\documentclass[12pt,reqno]{amsart}
\usepackage{amsmath, amssymb, amsthm}
\usepackage{url}
\usepackage[breaklinks]{hyperref}

\renewcommand{\Im}{\operatorname{Im}}

\renewcommand{\Re}{\operatorname{Re}}
\renewcommand{\Im}{\operatorname{Im}}

\renewcommand{\(}{\left\(}
\renewcommand{\)}{\right\)}
\renewcommand{\[}{\left\[}
\renewcommand{\]}{\right\]}

\numberwithin{equation}{section}
 \theoremstyle{plain}
\newtheorem{theorem}{Theorem}[section]
\newtheorem{lemma}[theorem]{Lemma}

\newtheorem{definition}[theorem]{Definition}

\newtheorem{corollary}[theorem]{Corollary}

   \makeatletter
\def\proof{\@ifnextchar[{\@oproof}{\@nproof}}
\def\@oproof[#1][#2]{\trivlist\item[\hskip\labelsep\textit{#2 Proof of\
#1.}~]\ignorespaces}
\def\@nproof{\trivlist\item[\hskip\labelsep\textit{Proof.}~]\ignorespaces}

\makeatother

\begin{document}
\title[A Lambert series associated with Siegel cusp forms]{An asymptotic expansion for a  Lambert series associated to Siegel cusp forms} 

\author{Babita}
\address{Babita\\Department of Mathematical Sciences\\ Indian Institute of Technology(Banaras Hindu University), Varanasi \\
221005, Uttar Pradesh, India.} 
\email{babita.rs.mat19@itbhu.ac.in}

\author{Abhash Kumar Jha}
\address{Abhash Kumar Jha\\Department of Mathematical Sciences\\Indian Institute of Technology(Banaras Hindu University), Varanasi \\
221005, Uttar Pradesh, India.} 
\email{abhash.mat@iitbhu.ac.in}

\author{Abhishek Juyal}
\address{Abhishek Juyal\\
Department of Mathematics,  \\
HNB Garhwal University,  BCC Campus,\\ 
Srinagar,  Uttarakhand, 246174.} 
\email{abhinfo1402@gmail.com}

 \author{Bibekananda Maji}
\address{Bibekananda Maji\\ Department of Mathematics \\
Indian Institute of Technology Indore \\
Indore, Simrol, Madhya Pradesh 453552, India.} 
\email{bibek10iitb@gmail.com, bibekanandamaji@iiti.ac.in}


\maketitle

\begin{abstract}
In 2000,  Hafner and Stopple proved  a conjecture of Zagier which states that the constant term of the automorphic function $|\Delta(x+iy)|^2$ i.e.,   the Lambert series $\sum_{n=1}^\infty \tau(n)^2 e^{-4 \pi n y}$ 
can be expressed in terms of the non-trivial zeros of the Riemann zeta function.  In this article,  we study a certain Lambert series associated to Siegel cusp forms and observe a similar phenomenon.  

\end{abstract}

\section{Introduction}
Let $\zeta(s)$  and  $\mu(n)$ denote  the Riemann zeta function and  the M\"{o}bius function, respectively. It is well-known that the reciprocal of the Riemann zeta function is the Dirichlet series associated to $\mu(n).$ Ramanujan,  in his second notebook \cite[p.~312]{Rama_2nd_Notebook},  \cite[eq. (37.3), p.~470]{BCB-V} mentioned the following identity associated to $\mu(n)$:

{\it For a positive real number $x$,  
\begin{align*}
\sum_{n=1}^{\infty} \frac{\mu(n)}{n} \exp\left({-\frac{x}{n^2}}\right) = \sqrt{\frac{\pi}{x}} \sum_{n=1}^{\infty} \frac{\mu(n)}{n} \exp\left(- \frac{\pi^2}{n^2 x } \right).  
\end{align*}}
However,  the above identity is not correct.  The corrected version of the above identity was obtained by Hardy and Littlewood \cite[p.~156]{HL-1916}.  Assuming all the non-trivial zeros of  $\zeta(s)$ are simple,  they proved that the following identity holds:
\begin{align}\label{Rama_Hardy_Little}
\sqrt{ \alpha} \sum_{n=1}^{\infty} \frac{\mu(n)}{n} \exp\left({- \frac{1}{\pi}\left(\frac{\alpha}{n}\right)^2}\right) 
& - \sqrt{\beta} \sum_{n=1}^{\infty} \frac{\mu(n)}{n} \exp\left({- \frac{1}{\pi}\left(\frac{\beta}{n}\right)^2}\right) \nonumber \\
& = -\frac{\sqrt{\pi}}{2\sqrt{\beta}} \sum_{\rho} \frac{ \Gamma\left(\frac{1-\rho}{2} \right) }{\zeta'(\rho)} \left( \frac{ \beta}{\sqrt{\pi}}\right)^\rho,
\end{align}
where  $\alpha$ and $\beta$ are two positive real numbers with $\alpha \beta = 1$ and the sum over $\rho$ runs through all  the non-trivial zeros of $\zeta(s)$.   The convergence of the series \eqref{Rama_Hardy_Little} is quite intricate.   In the same paper,  Hardy and Littlewood mentioned that they were unable to prove the convergence of this series even after assuming the Riemann hypothesis. However,  in  an indirect way,  they showed the convergence of this series  under the assumption of the following bracketing condition on the non-trivial zeros of $\zeta(s)$,  namely,
the non-trivial zeros $\rho_1$ and $\rho_2$ of $\zeta(s)$ are included in the same bracket if they satisfy
\begin{align}\label{bracketing condition}
|\Im(\rho_1) - \Im(\rho_2)| < \exp \left( -\frac{A_0 \Im(\rho_1)}{\log(\Im(\rho_1) )} \right) + \exp \left( -\frac{A_0 \Im(\rho_2)}{\log(\Im(\rho_2) )} \right),
\end{align}
where $A_0$ is some positive constant.  The identity  \eqref{Rama_Hardy_Little} also inspired Hardy and Littlewood to obtain the following equivalent criterion for the Riemann hypothesis (RH):
\begin{align}\label{Riesz type}
\sum_{n=1}^\infty  \frac{\mu(n)}{n} \exp\left({-\frac{x}{n^2}}\right) =  O\left( x^{-\frac{1}{4}+ \epsilon } \right), \quad \mathrm{as}\,\, x \rightarrow \infty,
\end{align}
Over the time,  the above equivalent criterion for RH has inspired mathematicians to find different generalization of \eqref{Rama_Hardy_Little} and \eqref{Riesz type} for $L$-functions associated to several objects.  For more details on this topic we refer to \cite{AGM22,  BK21,  DGV21,  DRZ,  DRZ-character,  GM23,  GV22}. We would like to mention here that, we also encounter an infinite series similar to the series \eqref{Rama_Hardy_Little} in the main result of this paper.


Let $\Delta(z)$ be the Ramanujan delta cusp form. The Fourier series expansion of $\Delta(z)$ is given by $$\Delta(z)=\sum_{n=1}^\infty \tau(n)\exp( 2 n \pi i z), ~~z \in \mathcal{H}.$$
 The Fourier coefficients $\tau(n)$ is known as Ramanujan tau function.  In 1981,  Zagier \cite[p.~417]{Zag}, \cite[p.~271]{Zag92} speculated  that the constant term of the automorphic function $\mathcal{F}(z):=y^{12}|\Delta(x+iy)|^2$,  that is, the following Lambert series 
\begin{equation}\label{Lambert series_Zagier}
c_0(y):=y^{12} \sum_{n=1}^{\infty} \tau^2(n) \exp({-4 \pi ny}), 
\end{equation}
has an oscillatory behaviour when $y \rightarrow 0^{+}$.
Moreover,  Zagier conjectured that its asymptotic expansion can be written in terms of the non-trivial zeros of $\zeta(s)$. 
Another interesting point he observed that one can numerically obtain non-trivial zeros of $\zeta(s)$ using the asymptotic expansion of \eqref{Lambert series_Zagier} and the Ramanujan tau function. 
Further, he claimed that the constant term $c_0(y)$ satisfies the following asymptotic expansion:
\begin{align*}
c_0(y) \sim \mathcal A + \sum_{\rho} y^{1 - \frac{\rho}{2}} A_{\rho},
\end{align*}
where the sum over $\rho$ runs through the non-trivial zeros of $\zeta(s)$,  $\mathcal{A}= \frac{3}{\pi} \langle \Delta , \Delta \rangle$, and the constant $A_{\rho}$ depends on $\rho$.  Assuming the Riemann hypothesis,  the above asymptotic expansion and the oscillatory behaviour of $c_0(y)$  have been proved by Hafner and Stopple \cite{HS}.   A similar phenomenon for Hecke eigenforms  for the full modular group as well as for the congruence subgroups have been observed in \cite{CKM} and \cite{CJKM}.  Subsequently,  an asymptotic expansion of an infinite series associated to a Hecke-Maass eigenforms has been studied by Banerjee and Chakraborty \cite{BC-19}.  Recently, this problem in the case of Hilbert modular forms has been examined by Agnihotri \cite{Agni}.  To know more about similar problems inspired from the above conjecture of Zagier,  we refer  to \cite{JMS-21,  JMS-22,  MNS23, MSS-22}.  
 
 Siegel modular forms are generalization of elliptic modular forms to higher dimension. They were first introduced by Siegel to study quadratic forms. A Siegel modular forms admits a Fourier-Jacobi expansion and the Fourier-Jacobi coefficients are Jacobi forms.  Motivated from the above mentioned conjecture of Zagier,  in this paper, we study certain Lambert series associated with the Fourier-Jacobi coefficients of Siegel modular forms. Interestingly, a similar phenomenon has been observed as in the case of elliptic modular forms. We use the analytic properties of certain Dirichlet series associated with Siegel modular forms (studied by Kohnen and Skoruppa \cite{k-s}), bounds of the Whittaker function, and a special evaluation of the Meijer $G$-function to prove our result.

 In the next section,  we provide necessary backgrounds which are essential to define the main problem of our article. 

\section{Preliminaries}
Let $\mathbb C$ and $\mathcal H$ denote the complex plane and complex upper half-plane, respectively. For two matrices $A$ and $B$ of appropriate sizes, we define $A[B]:={B}^tAB,$ where $B^t$ is the transpose of the matrix $B.$ We define $J:=\left(\begin{array}{cc} 0 & {I}\\{-I} & 0 \end{array}\right),$ where $I$ denotes the $2\times 2$ identity matrix. The Siegel upper half plane $ \mathcal{H}_2$ of degree $2$ is defined as 
$$
\mathcal{H}_2 := \{ Z  \in M_{2 \times 2}(\mathbb{C})~|~Z^t=Z,  \Im (Z)~ \text{is positive definite}  \}. 
$$
The full Siegel modular group  $\Gamma_2:=Sp_{4}(\mathbb Z)$  of genus $2$ is defined as 
\begin{equation*}
\Gamma_2:=\left\{M \in M_{4 \times 4}(\mathbb{Z}),J[M]=J\right\}.
\end{equation*} 
If we write the matrix $M\in \Gamma_2$ as a block decomposition $M=\begin{pmatrix}
A  & B\\
C  & D\\
\end{pmatrix}$, where  $ A,B,C,D \in M_{2 \times 2}(\mathbb{Z}),$ then we have
\begin{equation*}
\Gamma_2:=\left\{\begin{pmatrix}
A  & B\\
C  & D\\
\end{pmatrix} \bigg| A,B,C,D \in M_{2 \times 2}(\mathbb{Z}),A^tC = C^tA,  B^t D= D^t B,A^tD-B^tC=I_2\right\}.
\end{equation*}
The group $\Gamma_2$ acts on $\mathcal{H}_2$ as follows: 
$$
\begin{pmatrix}
A & B\\
C & D\\
\end{pmatrix} \cdot Z=(AZ+B)(CZ+D)^{-1}.
$$
\begin{definition}
A complex-valued holomorphic function $F: \mathcal{H}_2 \rightarrow \mathbb{C}$ is said to be a Siegel modular form of weight $k\in \mathbb N$ and degree two if it satisfies 
$$
F (M \cdot Z) = {\det} (CZ+D)^{k} F( Z),
$$
for all $M = \begin{pmatrix}
A  & B\\
C  & D\\
\end{pmatrix} \in \Gamma_2$,  $Z \in \mathcal{H}_2,$ and having a Fourier expansion of the form 
\begin{equation}\label{summation}
F(Z)= \sum_{\substack{T\ge 0}} A_F(T) e^{2\pi i (tr(TZ))}, 
\end{equation}
where the summation runs over positive semidefinite half-integral $2\times 2$ matrices $T=(t_{ij})$ ($ i.e., 2t_{ij},\;t_{ii}\in \mathbb{Z}$)  and $tr(\cdot)$ denotes the trace function.  

Further, we say  $F$ is a Siegel cusp form if and only if the summation in \eqref{summation} runs over positive definite half-integral matrices $T.$
\end{definition}
We denote the space of Siegel modular forms and Siegel cusp forms of weight $k$ and degree $2$ on $\Gamma_2$ by $M_k(\Gamma_2)$ and $S_k(\Gamma_2),$ respectively.
\noindent The Petersson scalar product on $S_k(\Gamma_2)$ is defined as
$$
\langle F, G \rangle:=\int_{\Gamma_2\setminus \mathcal H_2} F(Z)\overline{G(Z)}  (det~Y)^k dZ ,
$$
where $Z=X+iY$ and  $dZ=(det~Y)^{-3}dXdY$ is an invariant measure under the action of $\Gamma_2$ on $\mathcal H_2.$ The space $(S_k(\Gamma_2), \langle, \rangle)$ is a finite dimensional Hilbert space. For more details on the  theory on Siegel modular forms, we refer to \cite{1-2-3, maass}.

\subsection{Fourier-Jacobi expansion} Let $F(Z) \in M_k(\Gamma_2)$ be a Siegel modular form with the Fourier series expansion as in \eqref{summation}.

Write  $Z= \begin{pmatrix}
\tau & z\\
z & \tau'\\
\end{pmatrix} \in \mathcal{H}_2$ where $\tau, \tau'\in\mathcal H,~z\in \mathbb C.$ A matrix $T$ in the summation  \eqref{summation} can be written as $T=\begin{pmatrix}
m & r/2\\
r/2 & n\\
\end{pmatrix}$ with  $m,n \in \mathbb{N}, r \in \mathbb{Z}$ and   $4mn-r^2\ge 0.$
Then the Fourier series expansion  \eqref{summation} becomes
\begin{eqnarray*}\label{f-j}\nonumber
F(Z)= \sum_{T \ge 0} A_F(T) \exp({2\pi i(tr(TZ))})  & =& \sum\limits_{\substack{m,n,r\\4mn-r^2\ge 0}} A_F  \left( \begin{pmatrix}
m & r/2\\
r/2 & n\\
\end{pmatrix} \right)  \exp({2\pi i ( m \tau +r z + n \tau')})  \\
&=&  \sum\limits_{n=0}^{\infty} \phi_{n}(\tau,  z) \exp({2\pi i n \tau' }),
\end{eqnarray*}
where 
\begin{align*}
\phi_{n}(\tau,  z) = \sum_{\substack{m,r\\ 4mn-r^2\ge 0}} A_F  \left( \begin{pmatrix}
m & r/2\\
r/2 & n\\
\end{pmatrix} \right)  \exp({2\pi i ( m \tau +r z)}).  
\end{align*}
The coefficient $\phi_{n}(\tau,  z)$ is called the $n$th Fourier-Jacobi coefficients of $F$ and is a Jacobi form of weight $k$ and index $n$.  Jacobi forms are natural generalization of modular forms. Jacobi forms were first systematically studied  by Eichler and Zagier and they played a key role in the proof of Saito-Kurokawa conjecture. For more details on the theory of Jacobi forms, we refer to \cite{e-z}. It is easy to check that if $F\in S_k(\Gamma_2),$ then its Fourier-Jacobi expansion is given by $F(Z)=\sum\limits_{n=1}^{\infty} \phi_{n}(\tau,  z) \exp({2\pi i n \tau' })$ and each of the Fourier-Jacobi coefficients are Jacobi cusp forms.

\subsection{Dirichlet series associated to Siegel cusp forms} We now define  the Rankin-Selberg  type Dirichlet series associated to Siegel cusp forms $F_1, F_2\in S_k(\Gamma_2)$ with the Fourier-Jacobi expansion as follows:
\begin{equation}\label{fj-exp}
F_1(Z)=\sum\limits_{n=1}^{\infty} \phi_{n}(\tau,  z) \exp({2\pi i n \tau' }),~~~~~F_2(Z)=\sum\limits_{n=1}^{\infty} \psi_{n}(\tau,  z) \exp({2\pi i n \tau' })
\end{equation}
\begin{align}\label{Rankin-Selberg}
D_{F_1,  F_2}(s):= \zeta(2s-2k+4) \sum_{n=1}^\infty \frac{\langle \phi_{n},  \psi_{n}  \rangle}{n^s},
\end{align}
where $\langle \phi_{n},  \psi_{n}  \rangle$ denotes the Petersson scalar product (defined in \cite[p.~27,  eq.  (13)]{e-z}) of Jacobi cusp forms $\phi_n$ and $\psi_n.$
It is well-known \cite[p.~544,  Lemma 1]{k-s} that the Petersson scalar product of the Fourier-Jacobi coefficients $\phi_n$ and $\psi_n$ satisfies
$$ \langle \phi_n,\psi_n \rangle = O(n^k).$$
Further,  Kohnen \cite[p.~718]{Koh93},   \cite[p.~134]{Koh11} made the following Ramanujan-Petersson conjecture: for any $\epsilon>0$,  
\begin{align}\label{Ram-Pet}
\langle \phi_n,  \phi_n \rangle = O_{F_1,  \epsilon} \left( n^{k-1+\epsilon} \right).  
\end{align}
Kohnen and Sengupta \cite{KS17} proved the above conjecture for  Hecke eigenform $F\in S_k(\Gamma_2)$ which is a Saito-Kurokawa lift.  Recently,  Kumar and Paul \cite{KP21} proved that the above conjecture is true on average for arbitrary degree.  Assuming \eqref{Ram-Pet},  one can clearly see that the Dirichlet series $D_{F_1,F_2}(s)$ is absolutely convergent for $\Re(s)>k$.

In the next section we state the main results of our paper.

\section{Main Results}
As already mentioned in the introduction,  we are interested to study certain Lambert series associated to Siegel cusp forms.  More precisely, for $\alpha>0$,  we study the following Lambert series:
\begin{align}\label{Lambert_Siegel}
\sum_{n=1}^\infty \langle \phi_n,\psi_n \rangle \exp(-4 \pi n \alpha),
\end{align}
where $\phi_n$ and $\psi_n$ are defined as in \eqref{fj-exp}.
Before stating the main result of this article,  we define an arithmetic function $a_{F_1,  F_2}(n)$ that satisfy the following generating function:
\begin{align}\label{series1}
 \sum_{n=1}^\infty \frac{a_{F_1,  F_2}(n)}{n^s}= \frac{D_{F_1, F_2}(s)}{\zeta(2s+1-2k)},
\end{align}
which is absolutely convergent for $\Re(s) >k$.  Let $W_{\kappa,  \mu}(z)$  be the well-known { Whittakar} function which is one of the solutions of the following second order differential equation:
\begin{align*}
\frac{d^2w}{dz^2} + \left( \frac{1-4\mu^2+4 \kappa z -z^2}{4 z^2}  \right)w=0. 
\end{align*}
With the above notation and definitions in our hand,  we are ready to state the main result of this article. 

\begin{theorem}\label{main theorem}
Let $F_1$ and $F_2$ be two Siegel cusp forms of weight $k$ and degree $2$  with the Fourier-Jacobi expansion given  as in \eqref{fj-exp} and satisfying \eqref{Ram-Pet}.   Let $\alpha$ and $\beta$ be two positive real numbers such that $\alpha \beta =1$.  Under the assumption of the simplicity hypothesis of the non-trivial zeros of $\zeta(s)$,  we have
\begin{align}
\sum_{n=1}^\infty \langle \phi_n,\psi_n \rangle \exp(-4 \pi n \alpha) & =  
  \frac{ \beta^{2k-2}}{\pi^{3/2}} \sum_{n=1}^\infty \frac{ a_{F_1,  F_2}(n)}{ (4 \pi n \beta)^{\frac{k-1}{2}}}  \exp(-2\pi n \beta)  W_{\frac{k}{2}, \frac{k}{2}-1}(4\pi n \beta) \nonumber \\
& + R_{k} + \sum_{\rho} \frac{\Gamma\left(\frac{\rho}{2}+k-2 \right) D_{F_1, F_2}\left(\frac{\rho}{2}+k-2 \right) (4 \pi \alpha)^{2-k-\frac{\rho}{2} } }{\zeta'(\rho)},  \nonumber 
\end{align}
where the sum over $\rho$ runs through non-trivial zeros of $\zeta(s)$ and satisfy the bracketing condition \eqref{bracketing condition} and the term $R_k$ is given by
\begin{align}\label{Residue_k}
R_{k}=\begin{cases}
0, & \text { if }  \langle F_1,  F_2 \rangle = 0, \\
\frac{ 90 \langle F_1,  F_2 \rangle }{ \pi^2 \alpha^k},  & \text { if } \langle F_1,  F_2 \rangle \neq 0.
\end{cases}
\end{align}

\end{theorem}
As an application of the above theorem,  we have the following asymptotic expansion of the Lambert series defined in \eqref{Lambert_Siegel}. 
\begin{corollary}\label{application_asymptotic}
Let $F_1$ and $F_2$ be as in Theorem \ref{main theorem}. 
If $ \langle F_1,  F_2 \rangle \neq 0$,  then for $\alpha \rightarrow 0^{+} $,  we have
\begin{align}\label{asymptotic}
\alpha^k \sum_{n=1}^\infty \langle \phi_n,\psi_n \rangle \exp(-4 \pi n \alpha) \sim \frac{ 90 \langle F_1,  F_2 \rangle }{ \pi^2}. 
\end{align}
\end{corollary}

In the next section,  we collect a few results which will be essential for the proof of our main results.
\section{Some Well-known Results}
We complete the  Dirichlet series $D_{F_1,  F_2}(s)$ associated to the Siegel cusp forms $F_1$ and $F_2$ of weight $k$ and degree $2$ as follows:
$$
D^{*}_{F_1,  F_2}(s) := (2\pi)^{-2 s} \Gamma(s) \Gamma(s-k+2) D_{F_1,  F_2}(s).  
$$
The analytic properties of the above completed Dirichlet series is given below. 
\begin{theorem} \cite{k-s}\label{Functional equn}
The Dirichlet series  $D_{F_1,F_2}(s)$ can be analytically continued to an entire function if $\langle F_1,  F_2 \rangle=0$.  However,  if $\langle F_1,  F_2 \rangle \neq 0$,  then $D_{F_1,F_2}(s)$ has a simple pole at $s=k$ with residue $\frac{4^k \pi^{k+2}}{(k-1)!} \langle F_1,  F_2 \rangle$. 
The functional equation of $D_{F_1,  F_2}(s)$ is given by
\begin{equation}
D^{*}_{F_1,  F_2}(s)= D_{F_1,  F_2}^{*}(2k-2-s).
\end{equation} 
\end{theorem}
The above result will play a crucial role in the proof of the Theorem {\rm \ref{main theorem}}.  Next,  we define a well-known special function known as Meijer $G$-function.

Let $m,n,p,q$ be integers with $0\leq m \leq q$, $0\leq n \leq p$.  Let $a_1, \cdots, a_p$ and $b_1, \cdots, b_q$ be $p+q$ complex numbers such that $a_i - b_j \not\in \mathbb{N}$, for $i \in [1, n]$,  $ j \in [1,  m]$.  The Meijer $G$-function \cite[p.~415]{NIST} is defined by the following line integral:
\begin{align}\label{MeijerG}
	G_{p,q}^{\,m,n} \!\left(  \,\begin{matrix} a_1,\cdots , a_p \\ b_1, \cdots , b_q \end{matrix} \; \Big| z \right) = \frac{1}{2 \pi i} \int_L \frac{\prod_{j=1}^m \Gamma(b_j + w) \prod_{j=1}^n \Gamma(1 - a_j -w) z^{-w}  } {\prod_{j=m+1}^q \Gamma(1 - b_j - w) \prod_{j=n+1}^p \Gamma(a_j + w)}\mathrm{d}w.  
\end{align}
Here, we assume that the line of integration $L$ separates the poles of the factors $\Gamma(b_j+w)$ from the poles of the factors $\Gamma(1-a_j-w)$. The above integral converges absolutely if $2(m+n)> p+q$
 and $ |\arg(z)| < (2m+2n - p-q) \frac{\pi}{2}$. 
We now state the following special evaluation of the Meijer $G$-function.  
\begin{lemma}\cite[p.~58]{KT} \label{Special case of Meijer G}
 For $|\arg(z)|< \frac{\pi}{2}$,  we have
\begin{align}
G_{1,2}^{\,2,0} \!\left(  \,\begin{matrix} a_1 \\ b_1,  b_2 \end{matrix} \; \Big| z \right) = z^{\frac{b_1+b_2-1}{2}} e^{-z/2} W_{\kappa, \mu}(z),
\end{align}
where $W_{\kappa, \mu}(z)$ is the Whittaker function with the variables $\kappa= \frac{1}{2}(b_1+b_2-1)-a_1$ and $\mu=\frac{1}{2}(b_1 - b_2)$.
\end{lemma}
The next lemma tells us about the behaviour of the Whittaker function for large argument.
\begin{lemma} \cite[p.~69,  eq. (2.11.29)]{NIST}  \label{Asym_Whittaker}
As $z \rightarrow \infty$ with $|\arg(z)| < 3 \frac{\pi}{2}$,  one has
\begin{align*}
W_{\kappa, \mu}(z) \sim \exp(-z/2) z^{\kappa}.
\end{align*}
\end{lemma}

\begin{lemma}[Stirling's bound for $\Gamma(z)$]\cite[p.~151,~A.4.]{IK}   \label{Stirling}
	Let $z = \sigma + i {T}$ with $p \leq \sigma \leq q.$ Then  we have
	\begin{equation}\label{Stirling_equn}
		|\Gamma (\sigma + i {T})| \sim \sqrt{2\pi} | {T}|^{\sigma - 1/2} e^{-\frac{1}{2} \pi |{T}|},   \quad {\rm as} \quad |{T}|\rightarrow \infty.
	\end{equation}
\end{lemma}

\begin{lemma}\cite[p.~219]{Tit}\label{bound_for_zeta}
Let $T_n$ be a sequence of positive real numbers such that $|T_n - \Im(\rho)| > \exp\left(  - C \frac{\Im(\rho)}{\log(\Im(\rho))} \right)$ for any non-trivial zeros $ \rho$ of $\zeta(s)$,  where $C$ is some positive constant.  Then
\begin{align*}
\frac{1}{|\zeta(x+iT_n)|}< \exp(C_1 T_n),
\end{align*}
where $0<C_1< \pi/4$.  
\end{lemma}

In the next section,  we present a proof of the Theorem {\rm \ref{main theorem}}.

\section{Proof of main results}
\begin{proof}[Theorem {\rm \ref{main theorem}}][]

Using the inverse Mellin transform of the gamma function $\Gamma(s)$,  we see that the inverse Mellin transform of the function $\frac{\Gamma(s) D_{F_1,  F_2}(s) (4\pi)^{-s}}{\zeta(2s-2k+4)}$ is the infinite series $
\sum\limits_{n=1}^\infty \langle \phi_n,\psi_n \rangle \exp(-4 \pi n \alpha). 
$
Therefore,  for some large positive constant $c>k$,   we have
\begin{align}\label{right_vertical_integral}
\sum_{n=1}^\infty \langle \phi_n,\psi_n \rangle \exp(-4 \pi n \alpha)
  =\frac{1}{2\pi i}\int_{(c)}\frac{\Gamma(s)  D_{F_1,  F_2}(s) (4 \pi \alpha)^{-s}}{\zeta(2s-2k+4)}  \mathrm{d}s.
\end{align}
Here and throughout the article the symbol $(c)$ denotes the line integral $c- i \infty$ to $c+ i \infty$.  Now we recall from the Theorem \ref{Functional equn} that the Rankin-Selberg Dirchlet series $D_{F_1,  F_2}(s)$ can be extended to an entire function if $\langle F_1,  F_2 \rangle=0$.  However, if $\langle F_1,  F_2 \rangle \neq 0$,  then $D_{F_1,F_2}(s)$ has a simple pole at $s=k$ with residue $\frac{4^k \pi^{k+2}}{(k-1)!} \langle F_1,  F_2 \rangle$.  Moreover,  from the functional equation of $D_{F_1,  F_2}(s)$ it is clear that the poles of $\Gamma(s)$ are neutralized by $D_{F_1,  F_2}(s)$.  Nevertheless, the non-trivial of zeros of $\zeta(2s-2k+4)$ will give us infinitely many poles of the integrand function in the strip $k-2 < \Re(s) < k - \frac{3}{2}$,  and the trivial zeros of $\zeta(2s-2k+4)$  are at $k-n$ for $n\geq 3$.  Thus,  we construct a closed rectangular contour $\mathcal{C}$ consisting of the end points $c-i T,  c+i T,  c_1 + i T,  c_1 - iT$,  where $T$ is some large positive constant.  Here, we consider $c_1 \in (k-3,  k-2)$ so that the non-trivial zeros of $\zeta(2s-2k+4)$ lie inside  the contour $\mathcal{C}$ and the trivial zeros of $\zeta(2s-2k+4)$ lie outside the contour $\mathcal{C}$.  Now,  applying the Cauchy's residue theorem, we obtain 
\begin{equation}\label{CRT}
\frac{1}{2\pi i} \int_{\mathcal{C}} \frac{\Gamma{(s)}D_{F_1,F_2}(s)(4 \pi \alpha)^{-s}}{\zeta(2s-2k+4)}   \mathrm{d}s  = R_{k} + \sum_{ |\Im(\rho)| < T} R_{\rho}, 
\end{equation}
where the term $R_k$ denotes the residual term at $s=k$ and $R_{\rho}$ denotes the residual term at the non-trivial zero $\rho$ of $\zeta(s)$. 
The residual term $R_k$ is given by
\begin{align}\label{residue_k}
R_{k}=\begin{cases}
0, & \text { if }  \langle F_1,  F_2 \rangle = 0, \\
\frac{ 90 \langle F_1,  F_2 \rangle }{ \pi^2 \alpha^k},  & \text { if } \langle F_1,  F_2 \rangle \neq 0.
\end{cases}
\end{align}
Here, we have used the fact that $\zeta(4)=\pi^4/90$.  Now we assume the grand simplicity hypothesis,  which states that the non-trivial zeros of the Riemann zeta function $\zeta(s)$ are simple.  This gives us the following evaluation of $R_\rho$,
\begin{align}\label{residue_non-trivial}
R_{\rho}= \frac{\Gamma\left(\frac{\rho}{2}+k-2 \right) D_{F_1, F_2}\left(\frac{\rho}{2}+k-2 \right) (4 \pi \alpha)^{2-k-\frac{\rho}{2} } }{\zeta'(\rho)}. 
\end{align}
Letting $T \rightarrow \infty$ in \eqref{CRT} and utilizing Stirling's bound \eqref{Stirling_equn}  for $\Gamma(s)$ and Lemma \ref{bound_for_zeta},  one can show that both the horizontal integrals will vanish.  Finally, using \eqref{right_vertical_integral},  we obtain
\begin{align}\label{after CRT}
\sum_{n=1}^\infty \langle \phi_n,\psi_n \rangle \exp(-4 \pi n \alpha) 
  =\frac{1}{2\pi i}\int_{(c_1)}\frac{\Gamma(s)  D_{F_1,  F_2}(s) (4 \pi \alpha)^{-s}}{\zeta(2s-2k+4)}  \mathrm{d}s + R_k + \sum_{\rho} R_{\rho}, 
\end{align}
where the sum over $\rho$ runs through the non-trivial zeros of $\zeta(s)$.  Here, we emphasize that the sum over $\rho$ is an infinite sum as we know the existence of infinitely many zeros of $\zeta(s)$ from Hardy's theorem \cite[p.~257]{Tit}.  However,  the convergence of this series is quite delicate as we do not know much about the lower bound for $\zeta'(s)$.  Thus,  under the assumption of convergence of the sum in \eqref{after CRT} one can proceed further.  Now our main aim is to simplify the following left vertical integral:
\begin{align}\label{left_vertical_integral}
V_k(\alpha):= \frac{1}{2\pi i}\int_{(c_1)}\frac{\Gamma(s)  D_{F_1,  F_2}(s) (4 \pi \alpha)^{-s}}{\zeta(2s-2k+4)}  \mathrm{d}s.  
\end{align}
From the functional equation of $D_{F_1, F_2}(s)$,  i.e.,  Theorem \ref{Functional equn},  one obtains
\begin{align}\label{Rankin_functional}
 \Gamma(s) D_{F_1,  F_2}(s)=\frac{(2 \pi)^{4 s-4k+4} \Gamma(2 k-s-2) \Gamma(k-s) D_{F_1,  F_2}(2 k-2-s)}{\Gamma(s-k+2)}.
\end{align}
Also,  the symmetric form of the functional equation of the Riemann zeta function  implies that
\begin{align}\label{zeta_functional}
\zeta(2s-2k+4) = \frac{\pi^{2s-2k+4} \Gamma\left( k-s -\frac{3}{2}  \right) \zeta(2k-2s-3)}{\sqrt{\pi} \Gamma(s-k+2 )}.
\end{align}
Now substituting \eqref{Rankin_functional} and \eqref{zeta_functional} in \eqref{left_vertical_integral},  we obtain
\begin{align}
V_k(\alpha)= \frac{ \sqrt{\pi}}{16^{k-1} \pi^{2k}}   \frac{1}{2\pi i}\int_{(c_1)}  \frac{\Gamma(2 k-2-s) \Gamma(k-s) D_{F_1,  F_2}(2 k-2-s)}{ \Gamma\left( k-s -\frac{3}{2}  \right) \zeta(2k-2s -3)} (4 \pi \beta)^s \mathrm{d}s,
\end{align}
where $\alpha  \beta =1$.  
To simplify further,  we make a change of variable $w=2k-2-s$ and upon simplification, we obtain
\begin{align}\label{change variable}
V_k(\alpha)= \frac{ \beta^{2k-2}}{\pi^{3/2}}   \frac{1}{2\pi i}\int_{(d_1)} \frac{\Gamma(w) \Gamma(w-k+2) D_{F_1, F_2}(w) (4\pi \beta)^{-w}}{\Gamma\left( w-k+ \frac{1}{2}  \right) \zeta(2w+1-2k)} \mathrm{d}w.  
\end{align}
One can easily verify that $k < \Re(w)=d_1 <k+1$ as we have assumed $ k-3 < \Re(s)=c_1 < k-2$ .  This allows us to write the infinite series expansions for $D_{F_1, F_2}(w)$ and $1/\zeta(2w+1-2k)$ as both the series are absolutely and uniformly convergent in the region $\Re(w)>k$. 
Let us write 
\begin{align}\label{series}
\frac{D_{F_1, F_2}(w)}{\zeta(2w+1-2k)} = \sum_{n=1}^\infty \frac{a_{F_1,  F_2}(n)}{n^w},~~\Re(w)>k.
\end{align}
Now plugging the above series expansion \eqref{series} in \eqref{change variable} and interchanging the summation and integration,  we obtain
\begin{align}\label{Final_V_k}
V_k(\alpha) = \frac{ \beta^{2k-2}}{\pi^{3/2}} \sum_{n=1}^\infty a_{F_1,  F_2}(n) I_{k}(n,  \beta),
\end{align}
where
\begin{align}
 I_{k}(n,  \beta):= \frac{1}{2\pi i}\int_{(d_1)} \frac{\Gamma(w) \Gamma(w-k+2)  (4\pi n \beta)^{-w}}{\Gamma\left( w-k+ \frac{1}{2}  \right)} \mathrm{d}w.  
\end{align}
Now invoking the definition \eqref{MeijerG} of the Meijer $G$-function with $m=2,  n=0, p=1,  q=2$ and $a_1= \frac{1}{2} -k,  b_1=0, b_2= 2-k$,  we see that the above integral is the following Meijer $G$-function:
\begin{align}
I_{k}(n,  \beta)= G_{1,2}^{2,0}\left(\begin{array}{c}
\frac{1}{2}-k \\
0,2-k
\end{array} \Big| 4 \pi n \beta \right). 
\end{align}
One can easily verify that  the above Meijer $G$-function does satisfy all the convergence criteria.
Now employing Lemma \ref{Special case of Meijer G},  we obtain
\begin{align}\label{Final_I_k}
I_{k}(n,  \beta)=  (4 \pi n \beta)^{\frac{1-k}{2}} W_{k/2,  k/2-1}(4\pi n \beta) \exp(-2\pi n \beta).
\end{align}
Substituting \eqref{Final_I_k} in \eqref{Final_V_k},  the final expression of the left vertical integral becomes
\begin{align}\label{Final expression_V_k}
V_{k}(\alpha) & =  \frac{ \beta^{2k-2}}{\pi^{3/2}} \sum_{n=1}^\infty a_{F_1,  F_2}(n) (4 \pi n \beta)^{\frac{1-k}{2}} W_{\frac{k}{2},  \frac{k}{2}-1}(4\pi n \beta) \exp(-2\pi n \beta). 
\end{align}
Now applying the asymptotic expansion of the Whittaker function (Lemma \ref{Asym_Whittaker}),  one can easily verify the convergence of the above infinite series since the arithmetic function $a_{F_1,  F_2}(n)$ has polynomial growth and the Whittaker function decays exponentially.  Finally,  
considering the expression \eqref{Final expression_V_k} of the left vertical integral $V_{k}(\alpha)$  in \eqref{after CRT}, together with residual terms \eqref{residue_k} and \eqref{residue_non-trivial},  we complete the proof of Theorem \ref{main theorem}. 
\end{proof}

\begin{proof}[Corollary {\rm \ref{application_asymptotic}}][]
In Theorem \ref{main theorem},  we have started with the assumption on $\alpha$,  $\beta$ that $\alpha \beta =1$.  Therefore,  $\alpha \rightarrow 0^{+}$ implies that $\beta \rightarrow \infty$.  Making use of Lemma \ref{Asym_Whittaker} and multiplying by $\alpha^k$,  it is easy to see that 
\begin{align}
 & \frac{  \alpha^k \beta^{2k-2}}{\pi^{3/2}} \sum_{n=1}^\infty \frac{ a_{F_1,  F_2}(n)}{ (4 \pi n \beta)^{\frac{k-1}{2}}}  \exp(-2\pi n \beta)  W_{\frac{k}{2}, \frac{k}{2}-1}(4\pi n \beta) \nonumber \\
  & \ll \beta^{k-2} \sum_{n=1}^\infty \frac{ a_{F_1,  F_2}(n)}{ (4 \pi n \beta)^{\frac{k-1}{2}}}  \exp(-4\pi n \beta)  (4\pi n \beta)^{\frac{k}{2}}  \ll \frac{1}{\beta^M}, 
\end{align}
for some large positive number $M$.  This shows that the above infinite series vanishes as $\beta \rightarrow \infty$.  Again,  multiplying by $\alpha^k$ in the infinite residual term in Theorem \ref{main theorem},  which involves non-trivial zeros of $\zeta(s)$,  we arrive at 
\begin{align}\label{sum_non-trivial}
\sum_{\rho} \frac{\Gamma\left(\frac{\rho}{2}+k-2 \right) D_{F_1, F_2}\left(\frac{\rho}{2}+k-2 \right) (4 \pi)^{2-k-\frac{\rho}{2} }  \alpha^{2- \frac{\rho}{2}}}{\zeta'(\rho)}.
\end{align}
It is well-known that the non-trivial zeros of $\zeta(s)$ lie in the critical strip $0<\Re(s)<1$.  Therefore,  we have $3/2 < \Re\left(2- \frac{\rho}{2}\right) <2$,  and this indicates that the above infinite sum \eqref{sum_non-trivial} vanishes as $\alpha \rightarrow 0^{+}$.  Thus,  applying above two observations in Theorem \ref{main theorem} and looking at the residual term \eqref{Residue_k},  we complete the proof of \eqref{asymptotic}.  

\end{proof}

\begin{section}{Acknowledgement}
The second author is partially supported by SERB Start-up Research Grant (File No. SRG/2022/000487)  and MATRICS grant (File No. MTR/2022/000659).  The last author's research  supported by SERB MATRICS grant (File No. MTR/2022/000545).  Both the authors sincerely thank SERB for the support. 
\end{section}


\begin{thebibliography}{99}


\bibitem{AGM22} A.~ Agarwal, M.~Garg, and B.~ Maji,  \emph{Riesz-type criteria for the Riemann hypothesis},  Proc.  Amer.  Math.  Soc.,  {\bf 150} (2022), 5151--5163.

\bibitem{Agni}  R.~Agnihotri,  \emph{Lambert series associated to Hilbert modular form},   Int.  J.  Number Theory,  {\bf 18}  (2022),  1335--1349.  

\bibitem{Andrianov} A. ~N.~Andrianov,  Quadratic Forms and Hecke Operators.  Grundlehren der Math.  Wissenschaften 286,  Springer-Verlag (1987).  


\bibitem{BK21} S.~Banerjee and R.~Kumar,  \emph{Equivalent criterion for the grand Riemann hypothesis associated to Maass cusp forms},  submitted for publication,  2021. arXiv:2112.08143

\bibitem{BCB-V} B.~C.~ Berndt, Ramanujan’s Notebooks, Part V, Springer-Verlag, New York, 1998.

\bibitem{BC-19} S.~Banerjee and K.~Chakraborty,  \emph{Asymptotic behaviour of a Lambert series \'{a}  la Zagier: Maass case},  Ramanujan J.,  {\bf 48} (2019), 567--575.






\bibitem{boch} {S. ~B\"{o}cherer and  W. ~Kohnen},  \emph{Estimates for Fourier coefficients of Siegel cusp forms},  Math.~Ann.,~ {\bf 297} (1993), no. 3, 499--517.

\bibitem{1-2-3}
{J. ~H. ~Bruinier, G. ~van der Geer, G. ~Harder,  and D. ~Zagier,} The 1-2-3 of modular forms. Lectures from the Summer School on Modular Forms and their Applications held in Nordfjordeid, June 2004. Edited by Kristian Ranestad. Universitext. Springer-Verlag, Berlin, 2008. 

\bibitem{CKM} K.~Chakraborty, S.~Kanemitsu, and B.~Maji,  \emph{Modular-type relations associated to the Rankin-Selberg $L$-function},   Ramanujan J., {\bf 42} (2017), 285--299.


\bibitem{CJKM} K.~Chakraborty, A.~Juyal, S.~D.~Kumar, and B.~Maji,  \emph{An asymptotic expansion of a Lambert series associated to cusp forms},   Int. J. Number Theory, {\bf 14} (2018), 289--299.



\bibitem{DGV21} A. ~Dixit, S. ~Gupta,  and A. ~Vatwani,  \emph{A modular relation involving non-trivial zeros of the Dedekind zeta function, and the generalized Riemann hypothesis},  J. Math. Anal. Appl.,  {\bf 515} (2022), no. 2, 126435. 

\bibitem{DRZ} A.~Dixit, A.~Roy, and A.~Zaharescu,  \emph{Ramanujan-Hardy-Littlewood-Riesz phenomena for Hecke forms},   J. Math. Anal. Appl.,  {\bf 426} (2015),  594--611.

\bibitem{DRZ-character} A.~Dixit, A.~Roy, and A.~Zaharescu,  \emph{Riesz-type criteria and theta transformation analogues},  J.  Number Theory, {\bf 160} (2016), 385--408.

\bibitem{e-z}
M. Eichler and D. Zagier, The Theory of Jacobi Forms, Progr. in Math. vol. {\bf 55}, Birkhauser, Boston(1985).

\bibitem{GM23} M.~Garg and B.~ Maji,  \emph{Hardy-Littlewood-Riesz type equivalent criteria for the generalized Riemann hypothesis},  Monatsh.  Math.,  (2023).  https://doi.org/10.1007/s00605-023-01857-8


\bibitem{GV22}
S. ~Gupta and A. ~Vatwani,  \emph{Riesz type criteria for L-functions in the Selberg class},  accepted,  Canad. J. Math.  arXiv:2211.02954.  









%
%

\bibitem{HS} J.~Hafner and  J.~Stopple,  \emph{A heat kernel associated to Ramanujan's tau function,}  Ramanujan J.,  {\bf 4} (2000), 123--128.

\bibitem{HL-1916} G.~H.~ Hardy and J.~E.~Littlewood,  \emph{Contributions to the theory of the Riemann zeta-function and the theory of the distribution of primes},  Acta Math.,   {\bf 41} (1916), 119--196.


\bibitem{IK} H.~Iwaniec,  E.~Kowalski, Analytic number theory, American Mathematical Society Colloquium Publications, vol. 53, 2004. 









\bibitem{KT}
S.~Kanemitsu and H.~Tsukada,  Contributions to the theory of zeta-functions.  The modular relation supremacy.  Series on Number Theory and its Applications,  10. World Sci.,  Singapore,  2015. 

\bibitem{JMS-21} A.~Juyal,  B. ~Maji, and  S.~Satyanarayana,  \emph{An exact formula for a Lambert series associated to a cusp form and the M\"{o}bious function},  Ramanujan J.,  {\bf  57} (2022),  769--784. 

\bibitem{JMS-22} A.~Juyal,  B. ~Maji,  and S.~Satyanarayana,  
\emph{An asymptotic expansion for a Lambert series associated to the symmetric
square $L$-function,}  Int. J. Number Theory,  {\bf 19} (2022), 553--567.

\bibitem{k-s}
{W. ~Kohnen and N. ~P. ~Skoruppa},  \emph{A certain Dirichlet series attached to Siegel modular forms of degree two,}  Invent. Math.,  {\bf 95} (1989), no. 3,  541--558.


\bibitem{kohnen2}
W.~ Kohnen,  \emph{Estimates for Fourier coefficients of Siegel cusp forms of degree two. II.}  {Nagoya Math. J.,} {\bf  128} (1992), 171--176

\bibitem{Koh93}
W.~Kohnen,  \emph{Jacobi forms and Siegel modular forms: recent results and problems}, Enseign.  Math.  {\bf 39} (1993),   121--136. 


\bibitem{kohnen3}
{W. ~Kohnen,} \emph{Estimates for Fourier coefficients of Siegel cusp forms of degree two},  {Compositio Math.,} {\bf 87} (1993), no. 2, 231--240. 

\bibitem{KS17}
W.~Kohnen,  J.~Sengupta, \emph{Bounds for Fourier-Jacobi coefficients of Siegel cusp forms of degree two},   L-functions and automorphic forms, 159--169,  (eds J.~Bruinier and W. Kohnen; Springer, Berlin,  2017) 159--169.  



\bibitem{Koh11}
W.~Kohnen,  \emph{On the growth of the Petersson norms of Fourier-Jacobi coefficients of Siegel cusp forms},  Bul. Lond. Math.  Soc.  {\bf 43} (2011),  2025-2052.  

\bibitem{MNS23}
B.~ Maji,  P.~ Naskar, and S.~Sathyanarayana,
\emph{A series associated to Rankin-Selberg L-function and modified K-Bessel function}, submitted for publication,  2023.

\bibitem{MSS-22} B.~ Maji,  S.~ Sathyanarayana,  and B. ~R. ~Shankar,
\emph{An asymptotic expansion for a twisted Lambert series associated
to a cusp form and the Mobius function: level aspect},  Results Math., {\bf 77},  123 (2022). 



\bibitem{KP21}
B.~Kumar and  B.~Paul,  \emph{Ramanujan-Petersson conjecture Fourier-Jacobi coefficient for Siegel cups forms},  Bul.  Lond.  Math.  Soc.,  {\bf 53} (2021),  274--284. 

\bibitem{maass}
{H. Maass,} {Siegel's modular forms and Dirichlet series},
 Lecture Notes in Mathematics, {\bf 216}, Springer-Verlag, Berlin-New York, 1971.

\bibitem{NIST} F.~W.~J.~Olver, D.~W.~Lozier, R.~F.~Boisvert,  C.~W.~Clark, eds., NIST Handbook of Mathematical Functions, Cambridge University Press, Cambridge, 2010.


\bibitem{Rama_2nd_Notebook}
S.~Ramanujan,  Notebooks of Ramanujan,  Vol 2,  Tata Institute of FundamentaI Research,  Bombay, 1957.

 
 



 \bibitem{Tit} E.~C.~Titchmarsh, The Theory of the Riemann Zeta-function, Clarendon Press, Oxford, 1986.
 
 
 \bibitem{Zag} D.~Zagier,  \emph{The Rankin-Selberg method for automorphic functions which are not of rapid decay,}  J. Fac. Sci. Univ. Tokyo IA Math. {\bf 28} (1981), 415--437.

\bibitem{Zag92} D.~Zagier, Introduction to modular forms, From number theory to physics (Les Houches, 1989), Springer, 1992, 238--291.

\end{thebibliography}
\end{document}